\documentclass[12pt]{article}
\usepackage{amsmath,amssymb,amscd}
\title{Gorenstein Liaison and ACM Sheaves}
\author{Marta Casanellas
\\ Elena Drozd 
\\ Robin Hartshorne}
\date{}

\begin{document}

\maketitle

\begin{abstract}
We study Gorenstein liaison of codimension two subschemes of an 
arithmetically Gorenstein scheme X. Our main result is a criterion for two such 
subschemes to be in the same Gorenstein liaison class, in terms of the category 
of ACM sheaves on X. As a consequence we obtain a criterion for X to have the 
property that every codimension 2 arithmetically Cohen-Macaulay subscheme is in 
the Gorenstein liaison class of a complete intersection. Using these tools we 
prove that every arithmetically Gorenstein subscheme of $\mathbb{P}^n$ is in the 
Gorenstein liaison class of a complete intesection and we are able to 
characterize the Gorenstein liaison classes of curves on a nonsingular quadric 
threefold in $\mathbb{P}^4$. 
\end{abstract}


\setcounter{section}{-1}
\section{Introduction}
\label{sec0}

Liaison using complete intersection schemes has been widely studied since the
seminal paper of Peskine and Szpiro \cite{PS}.  More recently the notion of
liaison by arithmetically Gorenstein schemes, introduced by Schenzel \cite{S},
has assumed a prominent role in the study of subschemes of codimension at least
three in ${\mathbb P}^n$ \cite{KMMNP}.  The growing literature in this subject is
surveyed in the book \cite{M} and the recent article \cite{MN} of Migliore and
Nagel.

While most of the work on liaison has focused on subschemes of projective space
${\mathbb P}^n$, there are some results on other ambient varieties.  Rao remarks
already in \cite{R2} that the so-called ``Rao correspondence'' (see
$(2.16)$ below) holds on an arithmetically Gorenstein scheme, and Bolondi
and Migliore \cite{BM} establish the Lazarsfeld--Rao property on a smooth
arithmetically Gorenstein scheme.  Both of these results use complete
intersection linkage in codimension two.  Nagel \cite{Na} obtains similar results
from a slightly different point of view.  The recent paper \cite{RTLRP} puts the
Rao correspondence and the Lazarsfeld--Rao property in its most general form to
date:  in particular these results hold for complete intersection biliaison on an
integral arithmetically Cohen--Macaulay scheme.

The motivation for the current work came from trying to find a context more
general than the well-known case of curves in ${\mathbb P}^3$, yet more special
than the problem of Gorenstein liaison for curves in ${\mathbb P}^4$, which is
still wide open \cite{SEGL}.  So we decided to study curves in a nonsingular
quadric hypersurface in ${\mathbb P}^4$.  This project has expanded considerably
since its inception, and has led to the thesis of the second author \cite{D}, the
paper on Gorenstein biliaison of the first and third author \cite{GBAS}, and the
present article.

While the theory of biliaison can be satisfactorily developed on arithmetically
Cohen--Macaulay (ACM) schemes \cite{RTLRP}, \cite{GBAS}, the study of single
liaisons is more related to questions of duality, and requires Gorenstein
hypotheses on the ambient space.  So in this paper we study Gorenstein liaison of
codimension two subschemes of an arithmetically Gorenstein (AG) scheme $X$.  (See
\S1 for precise definitions).  Our main result $(5.1)$ is a criterion for
two such subschemes to be in the same $G$-liaison equivalence class, in terms of
the category of arithmetically Cohen--Macaulay (ACM) sheaves on $X$ (see
Definition $(2.6)$).  As a consequence we obtain a criterion for $X$ to have the
property that every codimension two ACM subscheme is in the Gorenstein liaison
class of a complete intersection (glicci).  In particular, this applies to the
nonsingular quadric hypersurface in ${\mathbb P}^4$ $(6.1)$ and solves the
problem we started with.

As a byproduct of the tools we develop for handling Gorenstein liaison, we are
able to prove that every arithmetically Gorenstein subscheme of ${\mathbb P}^n$
is glicci $(7.1)$, answering a natural question about Gorenstein liaison
\cite[p.~11]{KMMNP}, \cite{C}.

In \S1 we recall the definition and basic properties of liaison.  In \S2 we adapt
the ${\mathcal E}$- and ${\mathcal N}$-type resolutions used in \cite{MDP} and
elsewhere to our situation, and recall Rao's theorem.  \S3 studies the behavior
of these resolutions under $G$-liaison, and contains our main technical result
$(3.4)$, which allows one to simplify an ${\mathcal N}$-type resolution by a
$G$-liaison.  \S4 studies syzygies of ACM sheaves and introduces the notion of a
double-layered sheaf that is used in stating our main theorem.  \S5 has the main
theorem giving a criterion for $G$-liaison equivalence.  In the case of an
integral
$3$-dimensional AG scheme it also gives a criterion for the even $G$-liaison
class to be determined by the Rao module, as is the case in ${\mathbb P}^3$
\cite{R1}.  Then
\S6 has the applications to quadric hypersurfaces in ${\mathbb P}^4$ and
${\mathbb P}^5$, and
\S7 has the theorem that AG schemes in ${\mathbb P}^n$ are glicci.

\section{Liaison}
\label{sec1}

Liaison theory has its roots in the late nineteenth century.  The modern theory
of liaison begins with the paper of Peskine and Szpiro \cite{PS}, has been
developed by \cite{R1}, \cite{R2}, \cite{S}, \cite{MDP}, \cite{No}, \cite{Na},
\cite{KMMNP} and many others.  A good summary is in the book of Migliore
\cite{M}.  For convenience we recall the basic definitions and results that we
will use in this paper.  A scheme is {\em locally Cohen--Macaulay} (CM) if all of
its local rings are Cohen--Macaulay rings.  A closed subscheme $X$ of ${\mathbb
P}^n$ is {\em arithmetically Cohen--Macaulay} (ACM) if its homogeneous coordinate
ring $S(X) = k[x_0,\dots,x_n]/I_X$ is a Cohen--Macaulay ring.  A scheme is {\em
locally Gorenstein} if its local rings are all Gorenstein rings.  A closed
subscheme $X$ of ${\mathbb P}^n$ is {\em arithmetically Gorenstein} if its
homogeneous coordinate ring $S(X)$ is a Gorenstein ring.

\bigskip
\noindent
{\bf Definition 1.1.} 
Let $V_1,V_2,Y$ be subschemes of ${\mathbb P}_k^n$, all equidimensional without
embedded components, and of the same dimension.  We say $V_1$ {\em is linked to}
$V_2$ {\em by} $Y$, in symbols $V_1\ {\underset {Y}{\sim}}\ V_2$, if
\begin{itemize}
\item[1)] $V_1,V_2$ are contained in $Y$ and
\item[2)] ${\mathcal I}_{V_2,Y} \cong {\mathcal H}om({\mathcal O}_{V_1},{\mathcal
O}_Y)$ and ${\mathcal I}_{V_1,Y} \cong {\mathcal H}om({\mathcal
O}_{V_2},{\mathcal O}_Y)$.
\end{itemize}

\bigskip
Note the relation of linkage is symmetric by definition.  The following
proposition gives the existence and some properties of the linked scheme.

\bigskip
\noindent
{\bf Proposition 1.2.}  {\em Let $V_1 \subseteq Y$ be equidimensional subschemes
of ${\mathbb P}^n$, of the same dimension, without embedded components, and
assume that $Y$ is locally Gorenstein.  Define a subscheme $V_2$ of $Y$ by
setting ${\mathcal I}_{V_2,Y} = {\mathcal H}om({\mathcal O}_{V_1},{\mathcal
O}_Y)$ with its natural embedding in ${\mathcal O}_Y$.  Then}
\begin{itemize}
\item[a)] {\em $V_2$ is also equidimensional of the same dimension without
embedded components (unless $V_1 = Y$ in which case $V_2$ is empty) and $V_1$ and
$V_2$ are linked by $Y$.}
\item[b)] {\em $V_1$ is locally {\em CM} if and only if $V_2$ is locally {\em
CM}.}
\item[c)] {\em If $Y$ is {\em AG}, then $V_1$ is {\em ACM} if and only if $V_2$
is {\em ACM}.}
\end{itemize}

\bigskip
\noindent
{\em Proof.}  See \cite{S} or \cite[\S5.2]{M}.  The hypothesis $Y$ locally
Gorenstein is essential, to make the sheaves ${\mathcal I}_{V_1,Y},{\mathcal
I}_{V_2,Y}$ be reflexive sheaves on $Y$.

\bigskip
\noindent
{\bf Definition 1.3.}  If $V_1$ and $V_2$ are linked by $Y$ and $Y$ is a complete
intersection in ${\mathbb P}^n$, we speak of a CI-{\em linkage}.  The equivalence
relation generated by CI-linkages is called CI-{\em liaison}.  If a CI-liaison is
accomplished by an even number of CI-linkages, we speak of {\em even} CI-{\em
liaison}.

If $Y$ is an arithmetically Gorenstein (AG) scheme in ${\mathbb P}^n$, we speak
analogously of $G$-{\em linkage}, $G$-{\em liaison} and {\em even $G$-liaison}.

If the linkage $V_1\ {\underset {Y}{\sim}}\ V_2$ takes place among subschemes of
some projective scheme $X$ in ${\mathbb P}^n$, we speak of {\em linkage in} $X$.

Assume now that $X$ is an AG scheme in $X$.  A {\em complete intersection in} $X$
is an intersection $Y = H_1 \cap \dots \cap H_r$ of divisors $H_i$ in $X$
corresponding to invertible sheaves ${\mathcal O}_X(a_i)$, $i = 1,\dots,r$, of
codimension $r$ in $X$.  Such a scheme $Y$ is AG in ${\mathbb P}^n$, so we can
use it for liaison in $X$.  In this case, we speak of CI-{\em linkage}, CI-{\em
liaison}, {\em even} CI-{\em liaison in} $X$.  Note that unless $X$ itself is a
complete intersection in ${\mathbb P}^n$, these linkages may not be CI-linkages
in ${\mathbb P}^n$.

If $Y \subseteq X$ is AG in ${\mathbb P}^n$, then we speak of $G$-{\em linkage},
$G$-{\em liaison} and {\em even $G$-liaison in} $X$.  These are also $G$-linkages
in ${\mathbb P}^n$, by definition.

\bigskip
\noindent
{\bf Definition 1.4.}  We recall the notion of biliaison from \cite{GDB}.  If
$V_1$ and $V_2$ are divisors on an ACM scheme $S$ of one higher dimension, and
there is a linear equivalence 
$V_2
\sim V_1 + mH$ as divisors on $S$, where $H$ is the hyperplane section, for some
$m \in {\mathbb Z}$, we say $V_2$ is obtained by an {\em elementary biliaison}
from $V_1$.  If $S$ is an ACM scheme satisfying $G_0$, we speak of a $G$-{\em
biliaison}.  If $S$ is a complete intersection scheme in ${\mathbb P}^n$ (resp.\
complete intersection in $X$), then we speak of CI-{\em biliaison} (resp.\
CI-{\em biliaison in} $X$).  The corresponding equivalence relation generated by
the elementary biliaisons is called {\em biliaison}, {\em $G$-biliaison}, CI-{\em
biliaison}, or CI-{\em biliaison in} $X$.

\bigskip
\noindent
{\bf Proposition 1.5.} a) CI-{\em biliaison in ${\mathbb P}^n$ is the same
equivalence relation as even {\em CI}-liaison.}

b) {\em If $X$ is an {\em AG} scheme in ${\mathbb P}^n$, then {\em CI}-biliaison
in $X$ is the same equivalence relation as even {\em CI}-liaison in $X$.}

c) {\em Any $G$-biliaison is an even $G$-liaison.}

\bigskip
\noindent
{\em Proofs.} a) \cite[4.4]{GD}.

b) Use the same proof as \cite[4.4]{GD}, replacing ${\mathbb P}^n$ by $X$.

c) \cite[3.3]{GDB}.

\bigskip
\noindent
{\bf Remark 1.6.} It is unknown whether $G$-biliaison is the same as even
$G$-liaison in ${\mathbb P}^n$, but there are examples of non-singular AG schemes
$X$ on which $G$-biliaison and even $G$-liaison of codimension $2$ subschemes are
not equivalent \cite[1.1]{GBAS}.

\bigskip
\noindent
{\bf Remark 1.7.}  Note that our definition of $G$-liaison in $X$ is not the same
as the one used by Nagel \cite{Na}.  He requires a Gorenstein ideal to be
perfect, with the result that on an AG scheme $X$, his $G$-liaison for
codimension two subschemes on $X$ reduces to our CI-liaison on $X$.

\section{Resolutions of codimension two subschemes}
\label{sec2}

The ${\mathcal E}$- and ${\mathcal N}$-type resolutions and their behavior under
liaison were introduced in \cite{MDP} for curves in ${\mathbb P}^3$.  They have
been generalized in \cite{BM}, \cite{No}, \cite{Na}, \cite{HMDP1}, \cite{RTLRP}
to higher dimensions and other ambient schemes.  The two types of resolutions can
be defined separately with minimal hypotheses on arbitrary schemes, and Rao's
theorem can be proved for the ${\mathcal N}$-type resolution on a projective ACM
scheme \cite{RTLRP}, but the relation between the two types of resolution under
liaison appears clearly only for Gorenstein schemes.  We review for convenience
the definitions and results we need in this paper.

\bigskip
\noindent
{\bf Hypotheses 2.1.}  Throughout this section we let $X$ denote a projective
equidimensional scheme with a fixed very ample invertible sheaf ${\mathcal
O}_X(1)$.  A coherent sheaf ${\mathcal L}$ on $X$ is {\em dissoci\'e} if it is
isomorphic to a direct sum $\oplus {\mathcal O}_X(a_i)$ for various $a_i \in
{\mathbb Z}$.  We denote by $C$ a closed subscheme of $X$ of codimension two,
with no embedded components, and let ${\mathcal I}_C = {\mathcal I}_{C,X}$ be the
ideal sheaf of $C$ in $X$.  For any coherent sheaf ${\mathcal F}$ on
$X$, we write $H_*^i({\mathcal F}) = \oplus_{n \in {\mathbb Z}} H^i(X,{\mathcal
F}(n))$.

\bigskip
\noindent
{\bf Definition 2.2.}  An ${\mathcal E}$-{\em type resolution} of ${\mathcal
I}_C$ is an exact sequence
\[
0 \rightarrow {\mathcal E} \rightarrow {\mathcal L} \rightarrow {\mathcal I}_C
\rightarrow 0
\]
of coherent sheaves on $X$, with ${\mathcal L}$ dissoci\'e and $H_*^1({\mathcal
E}) = 0$.

\bigskip
\noindent
{\bf Proposition 2.3.} {\em With $X,C$ as in $(2.1)$ assuming furthermore
$H_*^1({\mathcal O}_X) = 0$, an ${\mathcal E}$-type resolution of $C$ exists.}

\bigskip
\noindent
{\em Proof.} Just let $I_C$ be the homogeneous ideal of $C$ in $S(X)$, the
projective coordinate ring of $X$, let $L$ be a free graded $S(X)$ module mapping
onto $I_C$, and let $E$ be the kernel:  $0 \rightarrow E \rightarrow L
\rightarrow I_C \rightarrow 0$.  Sheafifying gives the required ${\mathcal
E}$-type resolution.  Note that $H_*^1({\mathcal E}) = 0$ because
$H_*^0({\mathcal L}) \rightarrow H_*^0({\mathcal I}_C)$ is surjective, being $L
\rightarrow I_C$ by construction and $H_*^1({\mathcal L}) = 0$ since ${\mathcal
L}$ is dissoci\'e and $H_*^1({\mathcal O}_X) = 0$ by hypothesis.

\bigskip
\noindent
{\bf Definition 2.4.} An ${\mathcal N}$-{\em type resolution} of $C$ is an exact
sequence
\[
0 \rightarrow {\mathcal L} \rightarrow {\mathcal N} \rightarrow {\mathcal I}_C
\rightarrow 0
\]
with ${\mathcal L}$ dissoci\'e and ${\mathcal N}$ coherent satisfying
$H_*^1({\mathcal N}^{\vee}) = 0$ and ${\mathcal E}xt^1({\mathcal N},{\mathcal
O}_X) = 0$.  (This property of ${\mathcal N}$ was called ``extraverti'' in
\cite{HMDP1} and \cite{RTLRP}.)

\bigskip
\noindent
{\bf Proposition 2.5.} {\em With $X,C$ as in $(2.1)$, if $X$ satisfies Serre's
condition $S_2$ and $H_*^1({\mathcal O}_X) = 0$, then an ${\mathcal N}$-type
resolution of $C$ exists.}

\bigskip
\noindent
{\em Proof.} \cite[1.12]{RTLRP}.

\bigskip
\noindent
{\bf Definition 2.6.} A coherent sheaf ${\mathcal E}$ on a scheme $X$ is {\em
locally Cohen--Macaulay} (locally CM) if for every point $x \in X$, depth
${\mathcal E}_x = \dim {\mathcal O}_x$.  If $X$ is an ACM scheme, we say
${\mathcal E}$ is an ACM {\em sheaf} if in addition $H_*^i({\mathcal E}) = 0$ for
$0 < i < \dim X$ (cf. \cite{GBAS}).

\bigskip
\noindent
{\bf Proposition 2.7.}  {\em Suppose that $X$ is an {\em AG} scheme, and $C$ a
codimension two subscheme with ${\mathcal E}$- and ${\mathcal N}$-type resolutions
as defined above.  Then}
\begin{itemize}
\item[a)] {\em $C$ is locally {\em CM} if and only if ${\mathcal E}$
(resp.~${\mathcal N}$) is locally {\em CM}.}
\item[b)] {\em $C$ is {\em ACM} if and only if ${\mathcal E}$ (resp.~${\mathcal
N}$) is {\em ACM}.}
\end{itemize}

\bigskip
\noindent
{\em Proof.} a) Since $X$ is AG, it is locally CM, and so the dissoci\'e sheaves
are locally CM.  Now $C$ locally CM equivalent to ${\mathcal E}$ locally CM
follows by a chase of depth across the exact sequence $0 \rightarrow {\mathcal
I}_C \rightarrow {\mathcal O}_X \rightarrow {\mathcal O}_C \rightarrow 0$ and the
defining sequence for ${\mathcal E}$.  The implication ${\mathcal N}$ locally CM
$\Rightarrow C$ locally CM is similar.  In the other direction, if $C$ is locally
CM, we obtain depth ${\mathcal N} \ge n - 1$ at each closed point of $X$, where
$n= \dim X$.  It remains to show $H_x^{n-1}({\mathcal N}) = 0$ for each closed
point.  Since $X$ is locally Gorenstein, this follows by local duality and the
hypothesis ${\mathcal E}xt^1({\mathcal N},{\mathcal O}) = 0$.

b) By part a) we may assume $C$ locally CM and ${\mathcal E}$ and ${\mathcal N}$
also locally CM.  Then $C$ is ACM if and only if $H_*^i({\mathcal I}_C) = 0$ for
$0 < i \le \dim C = n-2$.  Chasing the cohomology sequences, this is equivalent
to $H_*^i({\mathcal E}) = 0$ for $2 \le i < n$ and to $H_*^i({\mathcal N}) = 0$
for $1 \le i < n-1$.  The missing requirements are $H_*^1({\mathcal E}) = 0$,
which is in the definition of ${\mathcal E}$-type resolution, and
$H_*^{n-1}({\mathcal N}) = 0$, which by Serre duality on $X$ (see $(2.8)$ below)
is equivalent to
$H_*^1({\mathcal N}^{\vee}) = 0$, in the definition of ${\mathcal N}$-type
resolution.

\bigskip
\noindent
{\bf Proposition 2.8.} {\em Let $X$ be an {\em AG} scheme.  Then}
\begin{itemize}
\item[a)] {\em Every locally {\em CM} sheaf ${\mathcal E}$ is reflexive.}
\item[b)] {\em The functor ${\mathcal E} \rightarrow {\mathcal E}^{\vee}$ is an
exact contravariant functor on the category of locally {\em CM} sheaves.}
\item[c)] {\em Serre duality for ${\mathcal E}$ locally {\em CM} says
$H_*^i({\mathcal E}^{\vee})$ dual to $H_*^{n-1}({\mathcal E})$.}
\item[d)] {\em ${\mathcal E}$ is {\em ACM} if and only if ${\mathcal E}^{\vee}$
is {\em ACM}.}
\end{itemize}

\bigskip
\noindent
{\em Proof.} \cite[2.3]{GBAS} noting that on an AG scheme, $\mbox{Hom}({\mathcal
E},\omega)$ is a twist of ${\mathcal E}^{\vee}$.

\bigskip
\noindent
{\bf Proposition 2.9.}  {\em Let $X$ be an {\em AG} scheme, and $C$ a locally
{\em CM} subscheme of codimension $2$.}
\begin{itemize}
\item[a)] {\em If $\dim X \ge 3$, $C$ is {\em subcanonical} (i.e., the canonical
sheaf $\omega_C$ is isomorphic to ${\mathcal O}_C(\ell)$ for some $\ell \in
{\mathbb Z}$) if and only if it has an ${\mathcal N}$-type resolution with
${\mathcal N}$ of rank $2$.}
\item[b)] {\em If $\dim X \ge 2$, $C$ is {\em AG} if and only if it has an
${\mathcal N}$-type resolution with ${\mathcal N}$ {\em ACM} of rank $2$.}
\end{itemize}

\bigskip
\noindent
{\em Proof.}  a) This is the usual ``Serre correspondence'' and the usual proof
(see, e.g., \cite[1.1]{SVB}) works also in our case.  If $C$ has a resolution
\[
0 \rightarrow {\mathcal L} \rightarrow {\mathcal N} \rightarrow {\mathcal I}_C
\rightarrow 0
\]
with ${\mathcal N}$ of rank $2$, then ${\mathcal L}$ has rank $1$, so ${\mathcal
L} = {\mathcal O}_X(-a)$ for some $a \in {\mathbb Z}$.  Taking ${\mathcal
H}om(\cdot,{\mathcal O}_X)$ we get
\[
0 \rightarrow {\mathcal O}_X \rightarrow {\mathcal N}^{\vee} \rightarrow
{\mathcal O}_X(a) \rightarrow {\mathcal E}xt_{{\mathcal O}_X}^1({\mathcal
I}_C,{\mathcal O}_X) \rightarrow 0.
\]
But ${\mathcal E}xt^1({\mathcal I}_C,{\mathcal O}_C) \cong {\mathcal
E}xt^2({\mathcal O}_C,{\mathcal O}_X) \cong \omega_C(-\ell)$, where $\ell$ is
such that $\omega_X = {\mathcal O}_X(\ell)$, using the hypothesis $X$ is AG. 
Hence $\omega_C(-\ell) \cong {\mathcal O}_C(a)$ and $C$ is subcanonical.

Conversely, if $C$ is subcanonical, reasoning backwards, ${\mathcal
E}xt^1({\mathcal I}_C,{\mathcal O}_X) \cong {\mathcal O}_C(a)$ for some $a$. 
Taking the global section $1 \in \mbox{Ext}^1({\mathcal I}_C,{\mathcal O}_X(-a))
= H^0({\mathcal E}xt^1({\mathcal I}_C,{\mathcal O}_X(-a)))$ gives an extension
\[
0 \rightarrow {\mathcal O}_X(-a) \rightarrow {\mathcal N} \rightarrow {\mathcal
I}_C \rightarrow 0.
\]
It follows, as in the proof of \cite[1.11]{RTLRP} that ${\mathcal N}$ satisfies
$H_*^1({\mathcal N}^{\vee}) = 0$ and ${\mathcal E}xt^1({\mathcal N},{\mathcal O})
= 0$, so this is an ${\mathcal N}$-type resolution with ${\mathcal N}$ of rank
$2$.

Note we used the hypothesis $\dim X \ge 3$ for the isomorphism ${\mathcal
E}xt^1({\mathcal I}_C,{\mathcal O}_X) \cong {\mathcal E}xt^2({\mathcal
O}_C,{\mathcal O}_X)$.

b) We apply the argument of part a) to the homogeneous coordinate ring $S(X)$,
which has dimension $\ge 3$ and to its quotient $S(C)$, and use the fact that a
graded ring is Gorenstein if and only if it is CM and the canonical module is
free.

\bigskip
To illustrate the previous proposition, we study when a codimension two AG
scheme
$C$ occurs as a divisor $mH-K$ on an ACM scheme $Y$ that is a divisor on $X$.  To
put this in a more general context, recall that if $Y$ is an ACM scheme
satisfying $G_0$, we can define the {\em anticanonical divisor} $M = M_Y$, given
by an embedding of $\omega_Y$ as a fractional ideal in the sheaf of total
quotient rings ${\mathcal K}_Y$, even if $Y$ does not have a well-defined
canonical divisor \cite[2.7]{GDB}.  Recall also that if $Y$ is an ACM scheme in
${\mathbb P}^n$ satisfying $G_0$, and if $C$ is an effective divisor on $Y$,
linearly equivalent to $mH+M_Y$, for some $m \in {\mathbb Z}$, then $C$ is AG in
${\mathbb P}^n$ \cite[3.4]{GDB}, \cite[5.2,5.4]{KMMNP}.

\bigskip
\noindent
{\bf Proposition 2.10.} {\em Let $X$ be an {\em AG} scheme, and let $C$ be an
{\em AG} subscheme of codimension $2$ in $X$.  Then the following conditions are
equivalent:}
\begin{itemize}
\item[(i)] {\em There is an {\em ACM} divisor $Y \subseteq X$ satisfying $G_0$
and containing $C$ and an integer $m$ so that $C \sim mH+M_Y$ on $Y$, where $M_Y$
is the anticanonical divisor.}
\item[(ii)] {\em $C$ has an ${\mathcal N}$-type resolution with ${\mathcal N}$ an
{\em ACM} sheaf of rank $2$ that is an extension of two rank $1$ {\em ACM}
sheaves on $X$.}
\end{itemize}

\bigskip
\noindent
{\em Proof.} (i) $\Rightarrow$ (ii).  Since $C \sim mH+M_Y$ on $Y$, we have
${\mathcal I}_{C,Y} \cong \omega_Y(-m)$.  On the other hand, comparing with the
ideal sheaf ${\mathcal I}_C$ of $C$ on $X$, we have an exact sequence
\begin{equation}
\label{eq1}
0 \rightarrow {\mathcal I}_Y \rightarrow {\mathcal I}_C \rightarrow {\mathcal
I}_{C,Y} \rightarrow 0.
\end{equation}
We combine this with the natural exact sequence
\begin{equation}
\label{eq2}
0 \rightarrow {\mathcal O}_X \rightarrow {\mathcal O}(Y) \rightarrow \omega_Y
\otimes \omega_X^{\vee} \rightarrow 0
\end{equation}
of \cite[2.10]{GD}.  (Here we write ${\mathcal O}(Y)$ for the notation ${\mathcal
L}(Y)$ of \cite{GD}.)  Since $X$ is an AG scheme, we can write $\omega_X \cong
{\mathcal O}_X(\ell)$ for some $\ell \in {\mathbb Z}$.  Twisting by $a = \ell -
m$ we get
\begin{equation}
\label{eq3}
0 \rightarrow {\mathcal O}_X(a) \rightarrow {\mathcal O}(Y+a) \rightarrow
\omega_Y(-m) \rightarrow 0.
\end{equation}
Since $\omega_Y(-m) \cong {\mathcal I}_{C,Y}$, we can do the fibered sum
construction with the sequence \eqref{eq1} above and obtain two short exact
sequences
\begin{equation}
\label{eq4}
0 \rightarrow {\mathcal O}_X(a) \rightarrow {\mathcal N} \rightarrow {\mathcal
I}_C \rightarrow 0,
\end{equation}
\begin{equation}
\label{eq5}
0 \rightarrow {\mathcal I}_Y \rightarrow {\mathcal N} \rightarrow {\mathcal
O}(Y+a) \rightarrow 0.
\end{equation}
The first is an ${\mathcal N}$-type resolution of ${\mathcal I}_C$, and the
second shows that ${\mathcal N}$ is an extension of two rank $1$ ACM sheaves on
$X$.

(ii) $\Rightarrow$ (i).  Conversely, suppose given an ${\mathcal N}$-type
resolution of the form \eqref{eq4} above, and suppose that ${\mathcal N}$ is an
extension
\begin{equation}
\label{eq6}
0 \rightarrow {\mathcal L} \rightarrow {\mathcal N} \rightarrow {\mathcal M}
\rightarrow 0
\end{equation}
where ${\mathcal L},{\mathcal M}$ are rank $1$ ACM sheaves on $X$.  The composed
map ${\mathcal L} \rightarrow {\mathcal I}_C \rightarrow {\mathcal O}_X$ shows
that ${\mathcal L}$ is isomorphic to the ideal sheaf ${\mathcal I}_Y$ of an ACM
divisor $Y$ containing $C$.  Then by comparing Chern classes we find ${\mathcal
M} \cong {\mathcal O}(Y+a)$.  Dividing the sequence \eqref{eq4} by ${\mathcal
I}_Y$ in the second and third place we obtain
\begin{equation}
\label{eq7}
0 \rightarrow {\mathcal O}_X(a) \rightarrow {\mathcal O}(Y+a) \rightarrow
{\mathcal I}_{C,Y} \rightarrow 0.
\end{equation}
Comparing with the sequence \eqref{eq2} above, we conclude that ${\mathcal
I}_{C,Y} \cong \omega_Y(a-\ell) = \omega_Y(-m)$.  Therefore $C \sim mH+M_Y$ on
$Y$.  Note from the isomorphism ${\mathcal I}_{C,Y} \cong \omega_Y(a-\ell)$ that
$\omega_Y$ is locally free at the generic points of $Y$, so $Y$ satisfies $G_0$.

\bigskip
\noindent
{\bf Remark 2.11.} We may ask, when does the extension \eqref{eq5} for ${\mathcal
N}$ in the above proof split?  It splits if and only if there is an effective
divisor $Z \sim -Y - aH$ such that $C$ is the scheme-theoretic intersection $Y
\cap Z$.  Indeed, if ${\mathcal N}$ splits, then the sequence
\begin{equation}
\label{eq8}
0 \rightarrow {\mathcal O}_X(a) \rightarrow {\mathcal O}(-Y) \oplus {\mathcal
O}(Y+a) \rightarrow {\mathcal I}_C \rightarrow 0
\end{equation}
gives a map ${\mathcal O}_X \rightarrow {\mathcal O}(-Y-a)$ defining an
effective divisor $Z$, and then the sequence \eqref{eq8} shows that $C = Y \cap
Z$.  Conversely, if there is such a $Z$, then the sequence \eqref{eq8} gives an
${\mathcal N}$-type resolution with ${\mathcal N}$ split.

\bigskip
\noindent
{\bf Example 2.12.} a) Let $P$ be a point on a nonsingular cubic surface $X$ in
${\mathbb P}^3$.  Then $P$ has an ${\mathcal N}$-type resolution on $X$
\[
0 \rightarrow {\mathcal O}_X \rightarrow {\mathcal N} \rightarrow {\mathcal I}_P
\rightarrow 0.
\]
One can show that for a general point $P \in X$, the sheaf ${\mathcal N}$ is not
an extension of ACM line bundles.  However, if $P$ lies on a line $L$, then $P
\sim -H-K_L$ on $L$, and ${\mathcal N}$ is an extension
\[
0 \rightarrow {\mathcal O}(-L) \rightarrow {\mathcal N} \rightarrow {\mathcal
O}(L) \rightarrow 0.
\]
But $-L$ is not effective, so the extension does not split.

b) Now take two points $Q,R$ on the cubic surface $X$.  They have an ${\mathcal
N}$-type resolution
\[
0 \rightarrow {\mathcal O}_X(-1) \rightarrow {\mathcal N} \rightarrow {\mathcal
I}_{Q+R} \rightarrow 0.
\]
Two general points $Q,R$ are contained in a twisted cubic curve $Y$, and then
$Q+R \sim -K_Y$.  Therefore ${\mathcal N}$ is an extension
\[
0 \rightarrow {\mathcal O}_X(-Y) \rightarrow {\mathcal N} \rightarrow {\mathcal
O}_X(Y-1) \rightarrow 0.
\]
But $H-Y$ is not effective, so this extension does not split.

On the other hand, if $Q,R$ lie on a conic $Y$, we get a similar sequence using
the conic $Y$.  In this case $H-Y$ is effective, equal to a line $L$.  Still,
for general $Q,R \in Y$, the sequence for ${\mathcal N}$ does not split.  Only
when $\{Q,R\} = Y \cap L$ does the sequence split.

\bigskip
\noindent
{\bf Proposition 2.13} (Behavior under liaison).  {\em Suppose that $X$ is {\em
AG}, and that $C$ has a resolution of the form}
\[
0 \rightarrow {\mathcal A} \rightarrow {\mathcal B} \rightarrow {\mathcal I}_C
\rightarrow 0
\]
{\em with ${\mathcal A}$ and $B$ locally {\em CM} sheaves on $X$ and
$H_*^1({\mathcal A}) = H_*^1({\mathcal B}^{\vee}) = 0$.  (For example, either an
${\mathcal E}$- or an ${\mathcal N}$-type resolution.)  Let $Y$ be a codimension
two complete intersection in $X$, containing $C$, with resolution}
\[
0 \rightarrow {\mathcal O}(-a-b) \rightarrow {\mathcal O}(-a) \oplus {\mathcal
O}(-b) \rightarrow {\mathcal I}_Y \rightarrow 0.
\]
{\em Let $C'$ be the curve linked to $C$ by $Y$ $(1.2)$.  Then $C'$ has a
resolution}
\[
0 \rightarrow {\mathcal B}^{\vee}(-a-6) \rightarrow {\mathcal A}^{\vee}(-a-b)
\oplus {\mathcal O}(-a) \oplus {\mathcal O}(-b) \rightarrow {\mathcal I}_{C'}
\rightarrow 0.
\]

\bigskip
\noindent
{\em Proof.} This is the usual mapping cone construction (see, e.g.,
\cite[p.~60]{MDP} or \cite[5.1.20]{M}).  Since $C \subseteq Y$ we have ${\mathcal
I}_Y \subseteq {\mathcal I}_C$.  Because of the hypothesis $H_*^1({\mathcal A}) =
0$, the defining sections of ${\mathcal I}_Y$ lift to ${\mathcal B}$ and we get a diagram
\[
\begin{array}{ccccccccc}
& & & & & & 0 \\
& & & & & & \uparrow \\
& & & & & & {\mathcal I}_{C,Y} \\
& & & & & & \uparrow \\
0 & \rightarrow & {\mathcal A} & \rightarrow & {\mathcal B} & \rightarrow & {\mathcal I}_C &
\rightarrow & 0 \\
& & \uparrow & & \uparrow & & \uparrow \\
0 & \rightarrow & {\mathcal O}(-a-b) & \rightarrow & {\mathcal O}(-a) \oplus
{\mathcal O}(-b) & \rightarrow & {\mathcal I}_Y & \rightarrow & 0 \\
& & & & & & \uparrow \\
& & & & & & 0
\end{array}
\]
The mapping cone gives
\[
0 \rightarrow {\mathcal O}(-a-b) \rightarrow {\mathcal A} \oplus {\mathcal O}(-a)
\oplus {\mathcal O}(-b) \rightarrow {\mathcal B} \rightarrow {\mathcal I}_{C,Y} \rightarrow
0.
\]
Now ${\mathcal I}_{C,Y} \cong {\mathcal H}om({\mathcal O}_{C'},{\mathcal O}_Y)$
by definition of linkage.  Since $X$ and $Y$ are AG, this sheaf is a twist
$\omega_{C'}(c)$ for some $c \in {\mathbb Z}$.  We apply the functor ${\mathcal
H}om(\cdot,{\mathcal O}_X)$ to the sequence above.  Since ${\mathcal A}$ and ${\mathcal B}$
are locally CM and $X$ is locally Gorenstein, these sheaves are acyclic for
${\mathcal E}xt_{{\mathcal O}_X}^i(\cdot,{\mathcal O}_X)$.  Thus we obtain
\[
0 \rightarrow {\mathcal B}^{\vee} \rightarrow {\mathcal A}^{\vee} \oplus {\mathcal O}(a)
\oplus {\mathcal O}(b) \rightarrow {\mathcal O}(a+b) \rightarrow {\mathcal
E}xt_{{\mathcal O}_X}^2(\omega_{C'}(c),{\mathcal O}_X) \rightarrow 0.
\]
But ${\mathcal E}xt_{{\mathcal O}_X}^2(\cdot,{\mathcal O}_X)$ is a dualizing
functor for sheaves on $Y$, so we get ${\mathcal O}_{C'}(a+b)$ on the right. 
Splitting the sequence and twisting by $-a-b$ gives the result.

\bigskip
\noindent
{\bf Corollary 2.14.} {\em If $X$ is {\em AG} and $C,C'$ are linked by the
complete intersection $Y$ in $X$, then an ${\mathcal E}$-type (resp.~${\mathcal
N}$-type) resolution for $C$ gives an ${\mathcal N}$-type (resp.~${\mathcal
E}$-type) resolution for $C'$ by $(2.13)$.}

\bigskip
\noindent
{\em Proof.} Since $X$ is AG, the condition $H_*^1({\mathcal A}) = 0$ and
$H_*^1({\mathcal B}^{\vee}) = 0$ give $H_*^1({\mathcal B}^{\vee})$ and $H_*^1({\mathcal A}^{\vee\vee})
= 0$.  The last conclusion ${\mathcal E}xt^1({\mathcal A}^{\vee},{\mathcal O}) =
0$ follows by local duality from the fact that ${\mathcal A}$ and hence
${\mathcal A}^{\vee}$ $(2.8)$ are locally CM sheaves.

\bigskip
\noindent
{\bf Definition 2.15.}  We say sheaves ${\mathcal F}$ and ${\mathcal G}$ on $X$
are {\em stably equivalent} if there exist dissoci\'e sheaves ${\mathcal
L},{\mathcal M}$ such that ${\mathcal F} \oplus {\mathcal L} \cong {\mathcal G}
\oplus {\mathcal M}$.  We say ${\mathcal F}$ is {\em orientable} if ${\mathcal
F}$ is locally free of constant rank $r$ in codimension $1$, and there exists a
closed subset $Z \subseteq X$ of codimension $\ge 2$ and an integer $\ell$ such
that $\Lambda^r{\mathcal F} \cong {\mathcal O}_X(\ell)$ on $X-Z$.

\bigskip
\noindent
{\bf Theorem 2.16} (Rao).  {\em Let $X$ be an integral {\em AG} scheme of
dimension $\ge 2$.  Then the following sets are in one-to-one correspondence:}
\begin{itemize}
\item[(i)] {\em Even {\em CI}-liaison classes in $X$ of codimension two locally
{\em CM} subschemes $C \subseteq X$.}
\item[(ii)] {\em Orientable locally {\em CM} sheaves ${\mathcal E}$ with
$H_*^1({\mathcal E}) = 0$ on $X$, up to stable equivalence and shift.}
\item[(iii)] {\em Orientable locally {\em CM} sheaves ${\mathcal N}$ with
$H_*^1({\mathcal N}^{\vee}) = 0$ on $X$, up to stable equivalence and shift.}
\end{itemize}

{\em The bijections are accomplished by sending each $C$ in (i) to the sheaf
${\mathcal E}$ in an ${\mathcal E}$-type resolution for (ii), or to the sheaf
${\mathcal N}$ in an ${\mathcal N}$-type resolution for (iii).}

\bigskip
\noindent
{\em Proof.}  \cite{R2}, \cite{No}, \cite{Na}, \cite[6.2.5]{M}, \cite{RTLRP}. 
The equivalence of (i) and (iii) by sending $C$ to its ${\mathcal N}$-type
resolution is proved under more general hypotheses in \cite[2.4]{RTLRP}.  Note
that each psi-equivalence class mentioned there contains extraverti sheaves
unique up to stable equivalence \cite[1.12]{RTLRP}, and that these are the
sheaves that appear in an ${\mathcal N}$-type resolution.

The bijection of (i) and (ii) follows by using a single CI-liaison to relate a
subscheme $C$ to another subscheme $C'$, and the fact that ${\mathcal E}$- and
${\mathcal N}$-type resolutions are interchanged by CI-liaison and duality
$(2.13)$.

\bigskip
\noindent
{\bf Remark 2.17.}  It follows from the theorem of course that the sets (ii) and
(iii) are in bijective correspondence.  One bijection between these two sets is
given by the functor ${\mathcal F} \mapsto {\mathcal F}^{\vee}$, but that is not
the bijection coming from this theorem.  A consequence of $(4.2)$ below is that
the bijection in this theorem sends ${\mathcal N} \in$ (iii) to ${\mathcal
N}^{\sigma} \in$ (ii) and ${\mathcal E} \in$ (ii) to ${\mathcal
E}^{\vee\sigma\vee} \in$ (iii), up to stable equivalence.  (See \S3 below for the
definition of ${\mathcal N}^{\sigma}$.)

\section{Behavior under Gorenstein liaison}
\label{sec3}

We consider an AG scheme $X$ and its locally CM codimension two subschemes.  In
this section we investigate how the ${\mathcal E}$- and ${\mathcal N}$-type
resolutions behave under a Gorenstein liaison.  The mapping cone construction (as
in $(2.13)$) does not work for an ${\mathcal E}$-type resolution.  However, given
an ${\mathcal N}$-type resolution of $C$ and a $G$-liaison by an AG scheme $Y$ to
another scheme $C'$, we can obtain an ${\mathcal N}$-type resolution of $C'$ with
a more complicated sheaf in the middle.  On the other hand, we prove a key result
about how to simplify an ${\mathcal N}$-type resolution by $G$-liaison when
${\mathcal N}$ contains a rank $2$ ACM factor.

\bigskip
\noindent
{\bf Definition 3.1.}  For any locally CM sheaf ${\mathcal F}$ on the AG scheme
$X$ we define the {\em syzygy sheaf} ${\mathcal F}^{\sigma}$ to be the
sheafification of the kernel of a minimal free presentation of $F =
H_*^0({\mathcal F})$ over $S(X)$.  In other words, let $L \rightarrow F
\rightarrow 0$ be a minimal free cover, and sheafify to get
\[
0 \rightarrow {\mathcal F}^{\sigma} \rightarrow {\mathcal L} \rightarrow
{\mathcal F} \rightarrow 0
\]
with ${\mathcal L}$ dissoci\'e and $H_*^0({\mathcal L}) \rightarrow
H_*^0({\mathcal F})$ surjective.  Then ${\mathcal F}^{\sigma}$ is also locally
CM, and $H_*^1({\mathcal F}^{\sigma}) = 0$.

\bigskip
\noindent
{\bf Proposition 3.2.}  {\em Let $X$ be an {\em AG} scheme, and let $C$ be a
locally {\em CM} codimension two subscheme with ${\mathcal N}$-type resolution}
\[
0 \rightarrow {\mathcal L} \rightarrow {\mathcal N} \rightarrow {\mathcal I}_C
\rightarrow 0,
\]
{\em where ${\mathcal N}$ is locally {\em CM} with $H_*^1({\mathcal N}^{\vee}) =
0$ $(2.7)$.  Let $Y$ be a codimension $2$ {\em AG} subscheme containing $C$, with
resolution}
\[
0 \rightarrow {\mathcal O}_X(-a) \rightarrow {\mathcal E} \rightarrow {\mathcal
I}_Y \rightarrow 0
\]
{\em where ${\mathcal E}$ is a rank $2$ {\em ACM} sheaf (cf.~$(2.9)$).  Then the
subscheme $C'$ linked to $C$ by $Y$ has an ${\mathcal N}$-type resolution}
\[
0 \rightarrow {\mathcal M}^{\vee} \rightarrow {\mathcal G} \rightarrow {\mathcal
I}_{C'}(a) \rightarrow 0
\]
{\em with ${\mathcal M}$ dissoci\'e, and where ${\mathcal G}$ is an extension}
\[
0 \rightarrow {\mathcal E}^{\vee} \oplus {\mathcal L}^{\vee} \rightarrow
{\mathcal G} \rightarrow {\mathcal N}^{\sigma\vee} \rightarrow 0.
\]

\bigskip
\noindent
{\em Proof.}  We repeat the cone construction as in $(2.13)$.  We have ${\mathcal
I}_Y \subseteq {\mathcal I}_C$, and the induced map ${\mathcal E} \rightarrow
{\mathcal I}_C$ lifts to ${\mathcal N}$, because the $\mbox{Ext}^1$ term in the
sequence
\[
\mbox{Hom}({\mathcal E},{\mathcal N}) \rightarrow \mbox{Hom}({\mathcal
E},{\mathcal I}_C) \rightarrow \mbox{Ext}^1({\mathcal E},{\mathcal L})
\]
is zero since ${\mathcal L}$ is dissoci\'e and ${\mathcal E}$ is ACM $(2.8)$.  So
the mapping cone exists, and dualizing the resulting sequence as before, we get
\[
0 \rightarrow {\mathcal N}^{\vee} \rightarrow {\mathcal E}^{\vee} \oplus
{\mathcal L}^{\vee} \rightarrow {\mathcal I}_{C'}(a) \rightarrow 0.
\]
This is neither an ${\mathcal N}$-type nor an ${\mathcal E}$-type resolution.

But now consider the syzygy sequence for ${\mathcal N}$, where ${\mathcal M}$ is
a minimal dissoci\'e cover:
\[
0 \rightarrow {\mathcal N}^{\sigma} \rightarrow {\mathcal M} \rightarrow
{\mathcal N} \rightarrow 0.
\]
Dualizing this gives $(2.8)$
\[
0 \rightarrow {\mathcal N}^{\vee} \rightarrow {\mathcal M}^{\vee} \rightarrow
{\mathcal N}^{\sigma\vee} \rightarrow 0.
\]
Let ${\mathcal G}$ be the fibered sum of ${\mathcal E}^{\vee} \oplus {\mathcal
L}^{\vee}$ and ${\mathcal M}^{\vee}$ over ${\mathcal N}^{\vee}$.  Then we get
exact sequences
\[
0 \rightarrow {\mathcal M}^{\vee} \rightarrow {\mathcal G} \rightarrow {\mathcal
I}_{C'}(a) \rightarrow 0
\]
and
\[
0 \rightarrow {\mathcal E}^{\vee} \oplus {\mathcal L}^{\vee} \rightarrow
{\mathcal G} \rightarrow {\mathcal N}^{\sigma\vee} \rightarrow 0.
\]
Note that $H_*^1({\mathcal N}^{\sigma}) = 0$ by construction, and
$H_*^1({\mathcal E}) = 0$ since ${\mathcal E}$ is ACM.  It follows that
$H_*^1({\mathcal G}^{\vee}) = 0$, so that we do indeed get the required
${\mathcal N}$-type resolution.

\bigskip
\noindent
{\bf Remark 3.3.}  Thus we see that in general, performing a $G$-liaison makes
the sheaf appearing in the ${\mathcal N}$-type resolution more complex, by taking
a syzygy dual and an extension by a rank $2$ ACM sheaf.  Our next key result
shows that we can also reverse this process.

\bigskip
\noindent
{\bf Proposition 3.4.}  {\em Let $X$ be a normal {\em AG} scheme, and let $C$ be
a codimension two subscheme with an ${\mathcal N}$-type resolution}
\[
0 \rightarrow {\mathcal L} \rightarrow {\mathcal N} \rightarrow {\mathcal I}_C
\rightarrow 0.
\]
{\em  Suppose that ${\mathcal N}$ belongs to an exact sequence}
\[
0 \rightarrow {\mathcal E} \rightarrow {\mathcal N} \rightarrow {\mathcal N}'
\rightarrow 0
\]
{\em where ${\mathcal E}$ is an orientable rank $2$ {\em ACM} sheaf on $X$, and
${\mathcal N}'$ is another locally {\em CM} sheaf of rank $\ge 2$.  Then there
exists a subscheme $D$ in the same even $G$-liaison class as $C$, with a
resolution}
\[
0 \rightarrow {\mathcal L}' \rightarrow {\mathcal N}' \rightarrow {\mathcal
I}_D(a) \rightarrow 0
\]
{\em for some $a \in {\mathbb Z}$.  (If we assume furthermore $H_*^1({\mathcal
N}'{}^{\vee}) = 0$, this will be an ${\mathcal N}$-type resolution for $D$.)}

\bigskip
\noindent
{\em Proof.} Choose $b \gg 0$ so that ${\mathcal E}(b)$, ${\mathcal N}(b)$, and
hence also ${\mathcal N}'(b)$, are generated by global sections.  Let rank
${\mathcal N} = r+1$, so that rank ${\mathcal N}' = r-1$.  Choose $r-2$
sufficiently general sections of ${\mathcal N}'(b)$ so that the cokernel is the
twisted ideal sheaf of a codimension two subscheme $D$, giving
\[
0 \rightarrow {\mathcal O}(-b)^{r-2} \rightarrow {\mathcal N}' \rightarrow
{\mathcal I}_D(a) \rightarrow 0
\]
for some $a \in {\mathbb Z}$.  This is possible by repeated application of
\cite[2.6]{RTLRP}, since ${\mathcal N}'(b)$ is generated by global sections.

It follows that ${\mathcal I}_D(a+b)$ is also generated by global sections.  So
we can choose two of those sections defining a complete intersection scheme $Z$
in $X$, containing $D$.  Let $D'$ be the subscheme linked to $D$ by $Z$.  Forming
the diagram
\[
\begin{array}{ccccccccc}
0 & \rightarrow & {\mathcal O}(-b)^{r-2} & \rightarrow & {\mathcal N}' &
\rightarrow & {\mathcal I}_D(a) & \rightarrow & 0 \\
& & \uparrow & & \uparrow & & \uparrow \\
0 & \rightarrow & {\mathcal O}(-a-2b) & \rightarrow & {\mathcal O}(-b)^2 &
\rightarrow & {\mathcal I}_Z(a) & \rightarrow & 0
\end{array}
\]
as in $(2.13)$, the cone construction gives a resolution
\[
0 \rightarrow {\mathcal O}(-a-2b) \rightarrow {\mathcal O}(-b)^r \rightarrow
{\mathcal N}' \rightarrow {\mathcal R} \rightarrow 0
\]
where ${\mathcal R} = {\mathcal I}_{D,Z}(a)$.  Because of the linkage by $Z$,
${\mathcal R}$ is also equal to ${\mathcal H}om({\mathcal O}_{D'},{\mathcal
O}_Z)$, which is isomorphic to a twist $\omega_{D'}(d)$ for some $d$.  We split
this sequence in the middle by a sheaf ${\mathcal F}$, and write a diagram
\[
\begin{array}{ccccccccc}
& & 0 & & 0 \\
& & \uparrow & & \uparrow \\
0 & \rightarrow & {\mathcal F} & \rightarrow & {\mathcal N}' & \rightarrow &
{\mathcal R} & \rightarrow & 0 \\
& & \uparrow & & \uparrow \\
& & {\mathcal O}(-b)^r & \stackrel{\alpha}{\rightarrow} & {\mathcal N} \\
& & \uparrow & & \uparrow \\
& & {\mathcal O}(-a-2b) & \stackrel{\beta}{\rightarrow} & {\mathcal E} \\
& & \uparrow & & \uparrow \\
& & 0 & & 0
\end{array}
\]
whose left column and top row come from the resolution of ${\mathcal R}$; the
middle column is the given sequence with ${\mathcal N}$, the map $\alpha$ is
obtained by lifting the map of ${\mathcal O}(-b)^r \rightarrow {\mathcal N}'$ to
${\mathcal N}$, which is possible since ${\mathcal E}$ is an ACM sheaf; and
$\beta$ is obtained by restricting $\alpha$ to ${\mathcal O}(-a-2b)$.

Applying ${\mathcal H}om(\cdot,{\mathcal E})$ to the left-hand column, we obtain
\[
{\mathcal E}(b)^r \rightarrow {\mathcal E}(a+2b) \rightarrow {\mathcal
E}xt^1({\mathcal F},{\mathcal E}).
\]
Now $X$ is normal, so ${\mathcal F}$ is locally free in codimension $1$, so the
${\mathcal E}xt$ sheaf on the right has support in codimension $\ge 2$.  Let $W$
be the image of $H^0({\mathcal E}(b)^r)$ in $H^0({\mathcal E}(a+2b))$.  Then the
subsheaf ${\mathcal E}_0$ of ${\mathcal E}(a+2b)$ generated by $W$ is equal to
${\mathcal E}(d+2b)$ in codimension $1$.  Therefore \cite[2.6]{RTLRP} applies,
and a sufficiently general $s \in W$ will give a quotient ${\mathcal E}(a+2b)/s$
that is torsion-free and locally free in codimension $1$, hence the twisted ideal
sheaf of a subscheme of codimension $2$.  The proof of [loc. cit.] shows that for
$s \in W$ sufficiently general, the map $\beta+s$ will have the same effect.  Let
$t \in H^0({\mathcal E}(b)^r)$ be an element whose image is $s$.  Then we have a
new commutative diagram
\[
\begin{array}{ccccccccc}
& & 0 & & 0 & & 0 \\
& & \uparrow & & \uparrow & & \uparrow \\
0 & \rightarrow & {\mathcal F} & \rightarrow & {\mathcal N}' & \rightarrow &
{\mathcal R} & \rightarrow & 0 \\
& & \uparrow & & \uparrow & & \uparrow \\
0 & \rightarrow & {\mathcal O}(-b)^r & \stackrel{\alpha+t}{\rightarrow} &
{\mathcal N} & \rightarrow & {\mathcal G} & \rightarrow & 0 \\
& & \uparrow & & \uparrow & & \uparrow \\
0 & \rightarrow & {\mathcal O}(-a-2b) & \stackrel{\beta+s}{\rightarrow} &
{\mathcal E} & \rightarrow & {\mathcal I}_Y(c) & \rightarrow & 0 \\
& & \uparrow & & \uparrow & & \uparrow \\
& & 0 & & 0 & & 0
\end{array}.
\]
The left-hand and middle column and the top row are the same as before.  The
horizontal maps $\alpha+t$, $\beta+s$ are new, and ${\mathcal G}$ and ${\mathcal
I}_Y(c)$ are defined as their cokernels.  Now ${\mathcal I}_Y(c)$ is the twisted
ideal sheaf of a subscheme $Y$, by construction, since ${\mathcal E}$ is
orientable.  But ${\mathcal E}$ is also ACM of rank $2$, so $Y$ is an AG
subscheme of $X$.

The new sheaf ${\mathcal G}$ is locally free of rank $1$ in codimension $1$
because of the right-hand column and the fact that ${\mathcal R}$ has support in
codimension $2$.  On the other hand, from the middle row we see that depth
${\mathcal G} \ge 1$ at every point of codimension $2$.  Therefore ${\mathcal G}$
is torsion-free.  Thus ${\mathcal G}$ is the twisted ideal sheaf of a curve $C'$,
with the same twist, so ${\mathcal G} = {\mathcal I}_{C'}(c)$.

Now $C$ and $C'$ both have ${\mathcal N}$-type resolutions with the same sheaf
${\mathcal N}$, up to twist, so by Rao's theorem $(2.16)$, $C$ and $C'$ are in
the same even CI-liaison class, and a fortiori, in the same even $G$-liaison
class.  On the other hand $C'$ is contained in the AG scheme $Y$, and the diagram
we have just written shows that ${\mathcal R} = {\mathcal I}_{C',Y}(c)$. 
Therefore ${\mathcal R}$ also equals $\omega_{C''}(e)$, for some $e \in {\mathbb
Z}$, where $C''$ is the curve linked to $C'$ by $Y$.  But we have already seen
that ${\mathcal R} = \omega_{D'}(d)$, so $d = e$ and $C'' = D'$.  Thus we have
linkages $C'\ {\underset{Y}{\sim}}\ D'$ and $D'\ {\underset{Z}{\sim}}\ D$, and we
conclude that $C$ and $D$ are in the same even $G$-liaison class, as required.

\section{Syzygies and double-layered sheaves}
\label{sec4}

Throughout this section we consider locally CM sheaves on an AG scheme $X$.  We
have already defined $(3.1)$ the syzygy sheaf ${\mathcal F}^{\sigma}$ of a
locally CM sheaf ${\mathcal F}$.  In this section we develop some functorial
properties of this operation, and define the notion of double-layered sheaves,
which will be used in our main theorem.  The motivation for this curious
definition $(4.4)$ is that it comes from the ${\mathcal N}$-type resolution found
in
$(3.2)$ and it satisfies the property $(4.5)$ below.

\bigskip
\noindent
{\bf Proposition 4.1.}  {\em Let $X$ be an {\em AG} scheme.}
\begin{itemize}
\item[a)] {\em If ${\mathcal E}$ is locally {\em CM}, then ${\mathcal
E}^{\sigma}$ is locally {\em CM} and $H_*^1({\mathcal E}^{\sigma}) = 0$.}
\item[b)] {\em If $H_*^1({\mathcal E}^{\vee}) = 0$, then there is a dissoci\'e
sheaf ${\mathcal M}$ such that ${\mathcal E}^{\sigma\vee\sigma\vee} \oplus
{\mathcal M} \cong {\mathcal E}$.}
\item[c)] {\em ${\mathcal E}$ is dissoci\'e if and only if ${\mathcal E}^{\sigma}
= 0$.}
\item[d)] {\em If $0 \rightarrow {\mathcal E}' \rightarrow {\mathcal E}
\rightarrow {\mathcal E}'' \rightarrow 0$ is an exact sequence with
$H_*^1({\mathcal E}') = 0$, then there is a dissoci\'e sheaf ${\mathcal M}$ and
an exact sequence}
\[
0 \rightarrow {\mathcal E}'{}^{\sigma} \rightarrow {\mathcal E}^{\sigma} \oplus
{\mathcal M} \rightarrow {\mathcal E}''{}^{\sigma} \rightarrow 0.
\]
\end{itemize}

\bigskip
\noindent
{\em Proof.} These properties are all elementary and probably well-known, so we
just give the idea of proofs.
\begin{itemize}
\item[a)] ${\mathcal E}^{\sigma}$ is locally CM by chasing depth in the defining
sequence.   $H_*^1({\mathcal E}^{\sigma}) = 0$ since the map of associated modules
$L \rightarrow E \rightarrow 0$ in the definition is surjective, and $X$ is ACM.
\item[b)] Let
\[
0 \rightarrow {\mathcal E}^{\sigma} \rightarrow {\mathcal L} \rightarrow
{\mathcal E} \rightarrow 0
\]
be the defining sequence for ${\mathcal E}^{\sigma}$.  Taking duals gives
\cite[2.3]{GBAS}
\[
0 \rightarrow {\mathcal E}^{\vee} \rightarrow {\mathcal L}^{\vee} \rightarrow
{\mathcal E}^{\sigma\vee} \rightarrow 0.
\]
Since we have assumed $H_*^1({\mathcal E}^{\vee}) = 0$, the associated map of
modules is surjective on the right, though ${\mathcal L}^{\vee}$ is perhaps not
minimal.  Hence there is a dissoci\'e sheaf ${\mathcal M}$ such that ${\mathcal
E}^{\vee} \cong {\mathcal E}^{\sigma\vee\sigma} \oplus {\mathcal M}$.  Dualizing
gives the result (replacing ${\mathcal M}^{\vee}$ by ${\mathcal M}$).
\item[c)] ${\mathcal E}^{\sigma}$ is $0$ if and only if ${\mathcal E}$ is
isomorphic to the dissoci\'e sheaf ${\mathcal L}$ in the definition of ${\mathcal
E}^{\sigma}$.
\item[d)] Let ${\mathcal L}' \rightarrow {\mathcal E}'$ and ${\mathcal L}''
\rightarrow {\mathcal E}''$ be the minimal covers defining ${\mathcal
E}'{}^{\sigma}$ and ${\mathcal E}''{}^{\sigma}$.  The map ${\mathcal L}''
\rightarrow {\mathcal E}''$ lifts to ${\mathcal E}$ because of the hypotheses
$H_*^1({\mathcal E}') = 0$.  Thus ${\mathcal L}' \oplus {\mathcal L}''
\rightarrow {\mathcal E}$ is a map, surjective on the module level, so its kernel
is ${\mathcal E}^{\sigma}$ plus a dissoci\'e ${\mathcal M}$, which gives the
required exact sequence.
\end{itemize}

\bigskip
\noindent
{\bf Proposition 4.2.} {\em Let $X$ be an {\em AG} scheme of dimension $\ge 2$.}
\begin{itemize}
\item[a)] {\em If ${\mathcal E}$ is an {\em ACM} sheaf on $X$ $(2.6)$, then
${\mathcal E}^{\sigma}$ is also {\em ACM}.}
\item[b)] {\em If ${\mathcal E}$ is {\em ACM}, then ${\mathcal E}^{\sigma}$ has
no dissoci\'e direct summands.}
\item[c)] {\em If ${\mathcal E}$ is {\em ACM} and indecomposable, then ${\mathcal
E}^{\sigma}$ is also indecomposable.}
\end{itemize}

\bigskip
\noindent
{\em Proof.} a) The defining exact sequence for ${\mathcal E}^{\sigma}$ shows
that $H_*^i({\mathcal E}^{\sigma}) = 0$ for $2 \le i < \dim X$.  The case $i=1$
follows from $(4.1a)$.

b) Let $0 \rightarrow {\mathcal E}^{\sigma} \rightarrow {\mathcal L} \rightarrow
{\mathcal E} \rightarrow 0$ be the defining sequence for ${\mathcal E}^{\sigma}$,
and suppose ${\mathcal E}^{\sigma} = {\mathcal E}' \oplus {\mathcal M}$ with
${\mathcal M}$ dissoci\'e.  Then we can write
\[
\begin{array}{ccccccccc}
& & 0 & & 0 \\
& & \downarrow & & \downarrow \\
& & {\mathcal M} & = & {\mathcal M} \\
& & \downarrow & & \downarrow \\
0 & \rightarrow & {\mathcal E}^{\sigma} & \rightarrow & {\mathcal L} &
\rightarrow & {\mathcal E} & \rightarrow & 0 \\
& & \downarrow & & \downarrow & & \| \\
0 & \rightarrow & {\mathcal E}' & \rightarrow & {\mathcal F} & \rightarrow &
{\mathcal E} & \rightarrow & 0 \\
& & \downarrow & & \downarrow \\
& & 0 & & 0
\end{array}.
\]
Since ${\mathcal E}$ and ${\mathcal E}'$ are ACM, we find that $H_*^1({\mathcal
F}^{\vee}) = 0$.  Hence the middle column of this diagram splits, so ${\mathcal
F}$ is dissoci\'e and is a summand of ${\mathcal L}$.  But $H_*^1({\mathcal E}')
= 0$, so the bottom row contradicts minimality of ${\mathcal L}$.

c) If ${\mathcal E}^{\sigma}$ decomposes into ${\mathcal E}' \oplus {\mathcal
E}''$, then from $(4.1b)$ we get ${\mathcal E} \cong {\mathcal
E}'{}^{\vee\sigma\vee} \oplus {\mathcal E}''{}^{\vee\sigma\vee} \oplus {\mathcal
M}$.  But neither ${\mathcal E}'$ nor ${\mathcal E}''$ can be dissoci\'e, by b),
so ${\mathcal E}'{}^{\vee\sigma\vee}$ and ${\mathcal E}''{}^{\vee\sigma\vee}$ are
both non-zero $(4.1c)$, which shows that ${\mathcal E}$ is decomposable.

\bigskip
\noindent
{\bf Proposition 4.3.} {\em Let $X$ be {\em AG} as before, and let $C$ be a
locally {\em CM} codimension $2$ subscheme.}
\begin{itemize}
\item[a)] {\em If $C$ has an ${\mathcal N}$-type resolution}
\[
0 \rightarrow {\mathcal L} \rightarrow {\mathcal N} \rightarrow {\mathcal I}_C
\rightarrow 0,
\]
{\em then there is a dissoci\'e sheaf ${\mathcal M}$ and an ${\mathcal E}$-type
resolution}
\[
0 \rightarrow {\mathcal N}^{\sigma} \oplus {\mathcal L} \rightarrow {\mathcal M}
\rightarrow {\mathcal I}_C \rightarrow 0.
\]
\item[b)] {\em If $C$ has an ${\mathcal E}$-type resolution}
\[
0 \rightarrow {\mathcal E} \rightarrow {\mathcal L} \rightarrow {\mathcal I}_C
\rightarrow 0,
\]
{\em then there is a dissoci\'e sheaf ${\mathcal M}$ and an ${\mathcal N}$-type
resolution}
\[
0 \rightarrow {\mathcal M}^{\vee} \rightarrow {\mathcal E}^{\vee\sigma\vee}
\oplus {\mathcal L} \rightarrow {\mathcal I}_C \rightarrow 0.
\]
\end{itemize}

\bigskip
\noindent
{\em Proof.} a) Let
\[
0 \rightarrow {\mathcal N}^{\sigma} \rightarrow {\mathcal M} \rightarrow
{\mathcal N} \rightarrow 0
\]
be the defining sequence for ${\mathcal N}^{\sigma}$, with ${\mathcal M}$
dissoci\'e.  Then the kernel of the composed map ${\mathcal M} \rightarrow
{\mathcal I}_C$ is an extension of ${\mathcal L}$ by ${\mathcal N}^{\sigma}$,
which splits, because $H_*^1({\mathcal N}^{\sigma}) = 0$ $(4.1)$.

b) Let
\[
0 \rightarrow {\mathcal E}^{\vee\sigma} \rightarrow {\mathcal M} \rightarrow
{\mathcal E}^{\vee} \rightarrow 0
\]
be the defining sequence for ${\mathcal E}^{\vee\sigma}$.  Dualizing gives
\[
0 \rightarrow {\mathcal E} \rightarrow {\mathcal M}^{\vee} \rightarrow {\mathcal
E}^{\vee\sigma\vee} \rightarrow 0.
\]
If ${\mathcal G}$ is the fibered sum of ${\mathcal L}$ and ${\mathcal M}^{\vee}$
over ${\mathcal E}$, we get a sequence
\[
0 \rightarrow {\mathcal M}^{\vee} \rightarrow {\mathcal G} \rightarrow {\mathcal
I}_C \rightarrow 0
\]
where ${\mathcal G}$ is an extension
\[
0 \rightarrow {\mathcal L} \rightarrow {\mathcal G} \rightarrow {\mathcal
E}^{\vee\sigma\vee} \rightarrow 0.
\]
This extension splits because $H_*^1({\mathcal E}^{\vee\sigma}) = 0$ $(4.1)$, so
we get the ${\mathcal N}$-type resolution desired.

\bigskip
\noindent
{\bf Definition 4.4.}  A {\em double-layered} sheaf ${\mathcal E}$ on the AG
scheme $X$ is a sheaf ${\mathcal E}$ for which there exists a filtration
\[
0 = {\mathcal E}_0 \subseteq {\mathcal E}_1 \subseteq \dots \subseteq {\mathcal
E}_n = {\mathcal E}
\]
with the following two properties a),b):
\begin{itemize}
\item[a)] Each factor ${\mathcal E}_i/{\mathcal E}_{i-1}$ is either
\begin{itemize}
\item[(i)] a rank $2$ orientable ACM sheaf, not dissoci\'e, or
\item[(ii)] a sheaf ${\mathcal F}^{\sigma\vee}$, where ${\mathcal F}$ is of type
(i).
\end{itemize}
\item[b)] There exists an integer $r$ with $0 \le r \le n$, such that for $i \le
r$, ${\mathcal E}_i/{\mathcal E}_{i-1}$ is of type (i) and for $i > r$,
${\mathcal E}_i/{\mathcal E}_{i-1}$ is of type (ii).
\end{itemize}

Note that a double-layered sheaf is automatically ACM and orientable.

\bigskip
\noindent
{\bf Proposition 4.5.} {\em If ${\mathcal E}$ is a double-layered sheaf on the
{\em AG} scheme $X$, then there exists a dissoci\'e sheaf ${\mathcal M}$ such
that ${\mathcal E}^{\sigma\vee} \oplus {\mathcal M}$ is double-layered.}

\bigskip
\noindent
{\em Proof.} By induction on the least integer $n$ for which there exists a
filtration as in definition $(4.4)$.

\medskip
{\bf Case $n=1$.} In this case ${\mathcal E}$ itself is of type (i) or (ii) of
a).  If ${\mathcal E}$ is of type (i), then ${\mathcal E}^{\sigma\vee}$ is of
type (ii) by definition.  If ${\mathcal E}$ is of type (ii), then ${\mathcal E} =
{\mathcal F}^{\sigma\vee}$ for some sheaf ${\mathcal F}$ of type (i).  By $(4.1)$
there exists a dissoci\'e sheaf ${\mathcal M}$ such that ${\mathcal
E}^{\sigma\vee} \oplus {\mathcal M} = {\mathcal F}^{\sigma\vee\sigma\vee} \oplus
{\mathcal M} = {\mathcal F}$.  If ${\mathcal M} \ne 0$, then ${\mathcal F}$
having rank $2$ would be dissoci\'e, contradicting (i).  So ${\mathcal M} = 0$
and we find ${\mathcal E}^{\sigma\vee} = {\mathcal F}$ is of type (i).  Condition
b) of the definition is trivial, since there is only one factor in the
filtration.  So in this case ${\mathcal E}^{\sigma\vee}$ itself is double-layered.

\medskip
{\bf Case $n \ge 2$.}  We write
\[
0 \rightarrow {\mathcal E}_{n-1} \rightarrow {\mathcal E} \rightarrow {\mathcal
E}' \rightarrow 0
\]
where ${\mathcal E}' = {\mathcal E}_n/{\mathcal E}_{n-1}$.  Then ${\mathcal
E}_{n-1}$ and ${\mathcal E}'$ are both double-layered.  By the induction
hypothesis, there is a dissoci\'e sheaf ${\mathcal M}_1$ such that ${\mathcal
E}'' = {\mathcal E}_{n-1}^{\sigma\vee} \oplus {\mathcal M}_1$ is double-layered. 
On the other hand, ${\mathcal E}'$ has only a single factor, so by the Case $n=1$
${\mathcal E}'{}^{\sigma\vee}$ is itself double-layered with $n=1$.

Taking syzygies there is a dissoci\'e sheaf ${\mathcal M}_2$ and an exact
sequence $(4.1)$
\[
0 \rightarrow {\mathcal E}_{n-1}^{\sigma} \rightarrow {\mathcal E}^{\sigma}
\oplus {\mathcal M}_2 \rightarrow {\mathcal E}'{}^{\sigma} \rightarrow 0,
\]
whose dual gives
\[
0 \rightarrow {\mathcal E}'{}^{\sigma\vee} \rightarrow {\mathcal E}^{\sigma\vee}
\oplus {\mathcal M}_2^{\vee} \rightarrow {\mathcal E}_{n-1}^{\sigma\vee}
\rightarrow 0.
\]
Adding ${\mathcal M}_1$ to the middle and right-hand terms we get
\[
0 \rightarrow {\mathcal E}'{}^{\sigma\vee} \rightarrow {\mathcal E}^{\sigma\vee}
\oplus {\mathcal M}_2^{\vee} \oplus {\mathcal M}_1 \rightarrow {\mathcal E}''
\rightarrow 0.
\]

Now the filtration of ${\mathcal E}''$ together with this sequence gives a
filtration of the middle term ${\mathcal F} = {\mathcal E}^{\sigma\vee} \oplus
{\mathcal M}_2^{\vee} \oplus {\mathcal M}_1$ satisfying condition a) of the
definition of double-layered.  If ${\mathcal E}'$ was of type (ii), then
${\mathcal E}'{}^{\sigma\vee}$ is of type (i), and condition b) follows.  If
${\mathcal E}'$ was not of type (ii), then it is of type (i) and by condition b)
for
${\mathcal E}$, it follows that all the factors of ${\mathcal E}_{n-1}$ are also
of type (i).  But it is easy to see from Case $n=1$ and the induction so far that
the operation
${\mathcal E}$ goes to the double-layered sheaf ${\mathcal E}^{\sigma\vee} \oplus
{\mathcal M}$ constructed in this proof reverses the order and interchanges the
types (i) and (ii) of the factors.  Hence in this case all the factors of
${\mathcal E}''$ are of type (ii), and again condition b) is satisfied.  Thus
${\mathcal F}$ is double-layered, as required.

\bigskip
\noindent
{\bf Example 4.6.} If ${\mathcal F} \oplus {\mathcal M}$, where ${\mathcal F}$ is
ACM and ${\mathcal M}$ dissoci\'e, is isomorphic to a double-layered sheaf
${\mathcal E}$, the sheaf ${\mathcal E}$ and its factors are not uniquely
determined.  Here is a trivial example.  Take ${\mathcal F} = 0$.  This is double
layered with no factors.  On the other hand, let ${\mathcal E}$ be an orientable
rank $2$ ACM sheaf that is not dissoci\'e.  Then the syzygy sequence for
${\mathcal E}$, dualized, gives an exact sequence
\[
0 \rightarrow {\mathcal E}^{\vee} \rightarrow {\mathcal L}^{\vee} \rightarrow
{\mathcal E}^{\sigma\vee} \rightarrow 0.
\]
Thus ${\mathcal L}^{\vee}$ is double-layered with factors ${\mathcal E}^{\vee}$
and ${\mathcal E}^{\sigma\vee}$.  But ${\mathcal L}^{\vee}$ is of the form
${\mathcal F} \oplus {\mathcal M}$, taking ${\mathcal M} = {\mathcal L}^{\vee}$.

\section{The main theorem}
\label{sec5}

We begin with a criterion for two schemes to be in the same $G$-liaison class on
$X$.

\bigskip
\noindent
{\bf Proposition 5.1.} {\em Let $X$ be a normal {\em AG} scheme, and let $C$ be a
locally {\em CM} closed subscheme of codimension $2$ with an ${\mathcal N}$-type
resolution}
\[
0 \rightarrow {\mathcal L} \rightarrow {\mathcal N} \rightarrow {\mathcal I}_C
\rightarrow 0.
\]
{\em Write ${\mathcal N} = {\mathcal N}_0 \oplus {\mathcal M}_0$ with ${\mathcal
M}_0$ dissoci\'e and ${\mathcal N}_0$ having no dissoci\'e direct summands.  Then
another such subscheme $C'$ is in the same $G$-liaison class as $C$ on $X$ if and
only if it has an ${\mathcal N}$-type resolution}
\[
0 \rightarrow {\mathcal L}' \rightarrow {\mathcal N}' \rightarrow {\mathcal
I}_{C'}(a') \rightarrow 0
\]
{\em where ${\mathcal N}'$ satisfies the condition}
\begin{itemize}
\item[(*)] {\em There is an exact sequence}
\[
0 \rightarrow {\mathcal E}_1 \oplus {\mathcal L}_1 \rightarrow {\mathcal N}'
\rightarrow {\mathcal G} \rightarrow 0
\]
{\em where ${\mathcal L}_1$ is dissoci\'e, ${\mathcal E}_1$ is a rank $2$ {\em
ACM} sheaf, and ${\mathcal G}$ satisfies the condition}

\item[(**)] {\em There is a filtration}
\[
0 = {\mathcal G}_0 \subseteq {\mathcal G}_1 \subseteq \dots \subseteq {\mathcal
G}_n = {\mathcal G}
\]
{\em and an integer $1 \le r \le n$ such that}
\begin{itemize}
\item[1)] {\em for each $i<r$, the factor ${\mathcal G}_i/{\mathcal G}_{i-1}$ is
of type {\em (i)} of $(4.4)$, and}
\item[2)] {\em for $i=r$, the factor ${\mathcal G}_r/{\mathcal G}_{r-1}$ is
either ${\mathcal N}_0$ or ${\mathcal N}_0^{\sigma\vee}$, and}
\item[3)] {\em for $i>r$, the factor ${\mathcal G}_i/{\mathcal G}_{i-1}$ is of
type {\em (ii)} of $(4.4)$.}
\end{itemize}
\end{itemize}

\bigskip
\noindent
{\em Proof.} First suppose $C'$ is in the same $G$-liaison class as $C$.  The
proof is by induction on the number of $G$-liaisons needed to get from $C$ to
$C'$.  If this number is zero, we let ${\mathcal E}_1$ be a rank $2$ dissoci\'e
sheaf and use the ${\mathcal N}$-type resolution
\[
0 \rightarrow {\mathcal L} \oplus {\mathcal E}_1 \rightarrow {\mathcal N} \oplus
{\mathcal E}_1 \rightarrow {\mathcal I}_C \rightarrow 0.
\]
On the other hand, since ${\mathcal N} \cong {\mathcal M}_0 \oplus {\mathcal
M}_0$, we get a sequence
\[
0 \rightarrow {\mathcal E}_1 \oplus {\mathcal M}_0 \rightarrow {\mathcal N}
\oplus {\mathcal E}_1 \rightarrow {\mathcal N}_0 \rightarrow 0
\]
and it is clear that ${\mathcal N} \oplus {\mathcal E}_1$ satisfies the condition
(*).

For the induction step, suppose that $C'$ has a resolution
\[
0 \rightarrow {\mathcal L}' \rightarrow {\mathcal N}' \rightarrow {\mathcal
I}_{C'}(a') \rightarrow 0
\]
with ${\mathcal N}'$ satisfying (*), and suppose that $C''$ is obtained from $C'$
by one $G$-liaison, using a codimension $2$ AG scheme $Y$ with resolution
\[
0 \rightarrow {\mathcal O}_X(-a'') \rightarrow {\mathcal E} \rightarrow {\mathcal
I}_Y(a') \rightarrow 0.
\]
Then by $(3.2)$, $C''$ will have an ${\mathcal N}$-type resolution
\[
0 \rightarrow {\mathcal M}_1^{\vee} \rightarrow {\mathcal H} \rightarrow
{\mathcal I}_{C''}(a'') \rightarrow 0
\]
with ${\mathcal M}_1$ dissoci\'e and ${\mathcal H}$ belonging to an exact sequence
\[
0 \rightarrow {\mathcal E}^{\vee} \oplus {\mathcal L}'{}^{\vee} \rightarrow
{\mathcal H} \rightarrow {\mathcal N}'{}^{\sigma\vee} \rightarrow 0.
\]
At this point we need a lemma.

\bigskip
\noindent
{\bf Lemma 5.2.} {\em Suppose that ${\mathcal N}_0$ is a locally {\em CM} sheaf
with no dissoci\'e direct summands satisfying $H_*^1({\mathcal N}_0^{\vee}) =
0$, and suppose ${\mathcal N}'$ is a sheaf satisfying condition {\em (*)} of
$(5.1)$ with this ${\mathcal N}_0$.  Then there is a dissoci\'e sheaf ${\mathcal
M}$ such that
${\mathcal N}'{}^{\sigma\vee} \oplus {\mathcal M}$ satisfies condition {\em (**)}
of
$(5.1)$.}

\bigskip
\noindent
{\em Proof.}  Applying $(4.1)$ to the sequence
\[
0 \rightarrow {\mathcal E}_1 \oplus {\mathcal L}_1 \rightarrow {\mathcal N}'
\rightarrow {\mathcal G} \rightarrow 0
\]
and dualizing, we find there is a dissoci\'e sheaf ${\mathcal M}_1$ and an exact
sequence
\[
0 \rightarrow {\mathcal G}^{\sigma\vee} \rightarrow {\mathcal N}'{}^{\sigma\vee}
\oplus {\mathcal M}_1^{\vee} \rightarrow {\mathcal E}_1^{\sigma\vee} \rightarrow
0.
\]
On the other hand, the proof of $(4.5)$, with an extra factor ${\mathcal N}_0$ or
${\mathcal N}_0^{\sigma\vee}$ in the middle, shows that there is a dissoci\'e
sheaf ${\mathcal M}_1$ such that ${\mathcal G}' = {\mathcal G}^{\sigma\vee}
\oplus {\mathcal M}_2$ satisfies condition (**).  Note that since ${\mathcal
N}_0$ has no dissoci\'e direct summands, ${\mathcal N}_0^{\sigma\vee\sigma\vee} =
{\mathcal N}_0$ by $(4.2)$.  Thus we can write
\[
0 \rightarrow {\mathcal G}' \rightarrow {\mathcal N}'{}^{\sigma\vee} \oplus
{\mathcal M}_1^{\vee} \oplus {\mathcal M}_2 \rightarrow {\mathcal
E}_1^{\sigma\vee} \rightarrow 0,
\]
and then the filtration on ${\mathcal G}'$ together with ${\mathcal
E}_1^{\sigma\vee}$ gives a filtration on the middle term satisfying (**).

In the special case where ${\mathcal E}_1$ is dissoci\'e, the factor ${\mathcal
E}_1^{\sigma\vee}$ becomes $0$, and the proof is simpler.

\bigskip
\noindent
{\em Proof of $(5.1)$, continued.}  By the lemma, we can find a dissoci\'e sheaf
${\mathcal M}_2$ (new notation) such that ${\mathcal N}'{}^{\sigma\vee} \oplus
{\mathcal M}_2$ satisfies (**).  Then $C''$ has an ${\mathcal N}$-type resolution
\[
0 \rightarrow {\mathcal M}_1^{\vee} \oplus {\mathcal M}_2 \rightarrow {\mathcal
H} \otimes {\mathcal M}_2 \rightarrow {\mathcal I}_{C''}(a'') \rightarrow 0
\]
where the middle term ${\mathcal H} \oplus {\mathcal M}_2$ satisfies condition (*),
as required.

Conversely, suppose that $C'$ has a resolution of the form given in $(5.1)$.  We
proceed by induction on the number of non-dissoci\'e rank $2$ ACM sheaves
${\mathcal E}$ that appear as ${\mathcal E}_1$ in the exact sequence of
${\mathcal N}'$, or that appear as ${\mathcal E}$ or ${\mathcal E}^{\sigma\vee}$
in the factors of ${\mathcal G}$ for $i \ne r$.  If that number is zero, then
${\mathcal E}_1$ is dissoci\'e, and ${\mathcal G}$ is either ${\mathcal N}_0$ or
${\mathcal N}_0^{\sigma\vee}$.  The extension splits, since $H_*^1({\mathcal
N}_0^{\vee}) = 0$ and $H_*^1({\mathcal N}_0^{\sigma}) = 0$, so ${\mathcal N}'$ is
stably isomorphic to ${\mathcal N}_0$ or ${\mathcal N}_0^{\sigma\vee}$.  In the
first case, we conclude by Rao's theorem $(2.16)$ that $C'$ is in the even
CI-liaison class of $C_1$ and a fortiori in the same $G$-liaison class as $C$.  In
the second case, we do a single CI-liaison using a dissoci\'e rank $2$ sheaf
${\mathcal E}$ in $(3.2)$ to obtain a curve $C''$ whose ${\mathcal N}$-type
resolution involves a sheaf ${\mathcal N}''$ stably isomorphic to ${\mathcal
N}_0$, and then proceed as in the first case.

For the induction step we will use $(3.4)$.

\medskip
{\bf Case 1.} Suppose that ${\mathcal E}_1$ is not dissoci\'e.  Then we apply
$(3.4)$ with the sequence
\[
0 \rightarrow {\mathcal E}_1 \rightarrow {\mathcal N}' \rightarrow {\mathcal N}''
\rightarrow 0
\]
and obtain another scheme $D$ in the same $G$-liaison class as $C'$ having an
${\mathcal N}$-type resolution
\[
0 \rightarrow {\mathcal L}'' \rightarrow {\mathcal N}'' \rightarrow {\mathcal
I}_D(a) \rightarrow 0
\]
for some $a \in {\mathbb Z}$.  Now by construction ${\mathcal N}''$ belongs to an
exact sequence
\[
0 \rightarrow {\mathcal L}_1 \rightarrow {\mathcal N}'' \rightarrow {\mathcal G}
\rightarrow 0.
\]
Adding a rank $2$ dissoci\'e sheaf ${\mathcal E}'_1$ to ${\mathcal L}_1$ and to
${\mathcal N}''$, we find ${\mathcal N}'' \oplus {\mathcal E}'_1$ satisfies (*)
and has one fewer nondissoci\'e factor than $C'$, so the induction is
accomplished.

\medskip
{\bf Case 2.}  Suppose that ${\mathcal E}_1$ is dissoci\'e.  Then the sequence
for ${\mathcal N}'$ splits (since $H_*^1({\mathcal N}_0^{\vee}) = 0$, it being an
${\mathcal N}$-type resolution), so ${\mathcal N}' \cong {\mathcal E}_1 \oplus
{\mathcal L}_1 \oplus {\mathcal G}$.  If $r \ge 2$, we take the factor ${\mathcal
G}_1/{\mathcal G}_0$ as ${\mathcal E}$ in $(3.4)$ and proceed as in Case~1.  If
$r=1$, so that there is no factor of type (i) in ${\mathcal G}$, then we first
perform a single CI-liaison (as above), which converts type (ii) factors of
${\mathcal G}$ into type (i) factors, and then proceed as in the case $r>1$.

\bigskip
\noindent
{\bf Corollary 5.3.} {\em A codimension $2$ subscheme $C'$ of $X$ is in the
Gorenstein liaison class of a complete intersection (glicci) on $X$ if and only if
it has an ${\mathcal N}$-type resolution}
\[
0 \rightarrow {\mathcal L}' \rightarrow {\mathcal N}' \rightarrow {\mathcal
I}_{C'}(a') \rightarrow 0
\]
{\em where ${\mathcal N}' = {\mathcal N}_1 \oplus {\mathcal M}_1$ with ${\mathcal
N}_1$ a double-layered sheaf and ${\mathcal M}_1$ dissoci\'e.}

\bigskip
\noindent
{\em Proof.}  Taking $C$ to be a complete intersection in $X$ we have ${\mathcal
N}$ dissoci\'e, so ${\mathcal N}_0 = 0$.  Therefore the sheaf ${\mathcal G}$
appearing in the exact sequence
\[
0 \rightarrow {\mathcal E}_1 \oplus {\mathcal L}_1 \rightarrow {\mathcal N}'
\rightarrow {\mathcal G} \rightarrow 0
\]
of condition (*) is double-layered by condition (**).  If ${\mathcal E}_1$ is
dissoci\'e, the sequence splits and ${\mathcal N}' \cong {\mathcal G} \oplus
{\mathcal E}_1 \oplus {\mathcal L}_1$.  If ${\mathcal E}_1$ is not dissoci\'e, a
factor ${\mathcal L}_1$ splits off, so that ${\mathcal N}' \cong {\mathcal N}_1
\oplus {\mathcal L}_1$ and ${\mathcal N}_1$ belongs to a sequence
\[
0 \rightarrow {\mathcal E}_1 \rightarrow {\mathcal N}_1 \rightarrow {\mathcal G}
\rightarrow 0,
\]
which makes ${\mathcal N}_1$ double-layered, by definition.

\bigskip
\noindent
{\bf Theorem 5.4.} {\em Let $X$ be a normal {\em AG} scheme.  Then the following
conditions are equivalent:}
\begin{itemize}
\item[(D)] {\em Every codimension $2$ {\em ACM} subscheme $C$ is in the
Gorenstein liaison class of a complete intersection (glicci).}
\item[(E)] {\em Every orientable {\em ACM} sheaf on $X$ is stably equivalent to a
double-layered sheaf.}
\end{itemize}

\bigskip
\noindent
{\em Proof.}  (D) $\Rightarrow$ (E).  Suppose ${\mathcal N}$ is an orientable ACM
sheaf on $X$.  For $a \gg 0$ we can find a sequence
\[
0 \rightarrow {\mathcal O}(-a)^{r-1} \rightarrow {\mathcal N} \rightarrow
{\mathcal I}_C(b) \rightarrow 0
\]
where $r = \mbox{rank } {\mathcal N}$.  Then $C$ will be ACM $(2.7)$, and so by
condition (D) it will be glicci.  Then by $(5.3)$ it has an ${\mathcal N}$-type
resolution
\[
0 \rightarrow {\mathcal L}' \rightarrow {\mathcal N}' \rightarrow {\mathcal
I}_C(b) \rightarrow 0
\]
where ${\mathcal N}' = {\mathcal N}_1 \oplus {\mathcal M}_1$ with ${\mathcal
N}_1$ double-layered and ${\mathcal M}_1$ dissoci\'e.  But the ${\mathcal N}$ in
an ${\mathcal N}$-type resolution is uniquely determined up to stable equivalence
and shift by Rao's theorem $(2.16)$, so ${\mathcal N}$ is stably equivalent to a
double-layered sheaf.

(E) $\Rightarrow$ (D).  Given $C$ ACM of codimension $2$, let it have an
${\mathcal N}$-type resolution
\[
0 \rightarrow {\mathcal L} \rightarrow {\mathcal N} \rightarrow {\mathcal I}_C(a)
\rightarrow 0.
\]
Then ${\mathcal N}$ is orientable and ACM $(2.7)$, so by (E) it is stably
equivalent to a double-layered sheaf ${\mathcal N}_1$.  In other words, there are
dissoci\'e sheaves ${\mathcal M}_1$ and ${\mathcal M}_2$ such that ${\mathcal N}
\oplus {\mathcal M}_1 \cong {\mathcal N}_1 \oplus {\mathcal M}_2$.  Then the
${\mathcal N}$-type resolution
\[
0 \rightarrow {\mathcal L} \oplus {\mathcal M}_1 \rightarrow {\mathcal N}_1
\oplus {\mathcal M}_2 \rightarrow {\mathcal I}_C(a) \rightarrow 0
\]
for $C$ shows by $(5.3)$ that $C$ is glicci.

\bigskip
Next we specialize to the case of a $3$-dimensional normal AG scheme
$X$.  A {\em curve} will be a $1$-dimensional locally CM subscheme of $X$. 
Recall that for any curve $C$, we define the {\em Rao module} to be $M_C =
H_*^1({\mathcal I}_C)$.

\bigskip
\noindent
{\bf Theorem 5.5.} {\em Let $X$ be a normal $3$-dimensional {\em AG} scheme $X$,
and suppose that $X$ satisfies the (stronger) condition}

(E$'$) {\em Every orientable {\em ACM} sheaf on $X$ is stably equivalent to a
double-layered sheaf with all factors of type {\em (i)} in $(4.4)$.}

{\em Then two curves $C,C'$ on $X$ are in the same even $G$-liaison class if and
only if their Rao modules are isomorphic up to twist.}

\bigskip
\noindent
{\em Proof.}  One direction is well-known \cite[\S5.3]{M}.  For the other
direction, let $C$ be any curve on $X$, with Rao module $M$, and consider an
${\mathcal N}$-type resolution
\[
0 \rightarrow {\mathcal L} \rightarrow {\mathcal N} \rightarrow {\mathcal I}_C
\rightarrow 0
\]
of $C$.  Following the proof of \cite[4.7]{GBAS} we find another curve $C'$
depending only on $M$ and not on $C$, together with a resolution
\[
0 \rightarrow {\mathcal E} \rightarrow {\mathcal N} \rightarrow {\mathcal
I}_{C'}(a') \rightarrow 0
\]
where ${\mathcal E}$ is an orientable ACM sheaf.  By considering the syzygy sheaf
of ${\mathcal E}$, as in the proof of $(3.2)$, we obtain an ${\mathcal N}$-type
resolution of $C'$,
\[
0 \rightarrow {\mathcal M}^{\vee} \rightarrow {\mathcal G} \rightarrow {\mathcal
I}_{C'}(a') \rightarrow 0
\]
where ${\mathcal G}$ satisfies
\[
0 \rightarrow {\mathcal N} \rightarrow {\mathcal G} \rightarrow {\mathcal
E}^{\vee\sigma\vee} \rightarrow 0.
\]
Now ${\mathcal E}^{\vee}$ is an orientable ACM sheaf, so by hypothesis (E$'$) it
is stably equivalent to a double-layered sheaf ${\mathcal F}$ with all factors of
type (i).  Then by $(4.5)$ there is a dissoci\'e sheaf ${\mathcal M}_1$ such that
${\mathcal F}^{\sigma\vee} \oplus {\mathcal M}_1$ is double-layered with all
factors of type (ii).  Note also that ${\mathcal E}^{\vee\sigma\vee} = {\mathcal
F}^{\sigma\vee}$, since the operation $\sigma\vee$ kills off dissoci\'e factors.

On the other hand, let us write ${\mathcal N} = {\mathcal N}_0 \oplus {\mathcal
M}_0$ where ${\mathcal M}_0$ is dissoci\'e and ${\mathcal N}_0$ has no dissoci\'e
direct summands.  Then ${\mathcal G} = {\mathcal G}' \oplus {\mathcal M}_0$ where
${\mathcal G}'$ satisfies
\[
0 \rightarrow {\mathcal N}_0 \rightarrow {\mathcal G}' \oplus {\mathcal M}_1
\rightarrow {\mathcal F}^{\sigma\vee} \oplus {\mathcal M}_1 \rightarrow 0
\]
and hence ${\mathcal G}' \oplus {\mathcal M}_1$ satisfies condition (**) of
$(5.1)$.  Now ${\mathcal G}$ is stably equivalent to ${\mathcal G}' \oplus
{\mathcal M}_1$, and therefore $C'$ is in the same $G$-liaison as $C$ by
$(5.1)$.  But $C'$ depends only on the Rao module $M$, hence all curves with the
same Rao module (up to twist) are equivalent for $G$-liaison.  Since they have
the same Rao module (not their duals) it must be in the same even $G$-liaison
class.

\bigskip
\noindent
{\bf Remark 5.6.}  Of course it follows that two curves are connected by an odd
number of $G$-liaisons if and only if their Rao modules are dual (up to twist)
because a single $G$-liaison replaces the Rao module by its dual \cite[5.3.1]{M}.

\section{Quadric hypersurfaces}
\label{sec6}

\noindent
{\bf Theorem 6.1.} {\em If $X$ is a nonsingular quadric hypersurface of dimension
$3$ or $4$ in ${\mathbb P}^4$ or ${\mathbb P}^5$ respectively, then all
codimension $2$ {\em ACM} subschemes are glicci.}

\bigskip
\noindent
{\em Proof.} According to the theorems of Buchweitz, Eisenbud, Herzog, and
Kn\"orrer (see for example \cite[14.10]{Y}), $X$ has just one $(\dim X = 3)$ or
two $(\dim X = 4)$ indecomposable rank $2$ ACM sheaves, up to twist, and all
other ACM sheaves are direct sums of these and their twists and dissoci\'e
sheaves.  It follows that any ACM sheaf is stably equivalent to a double-layered
sheaf (in which all the extensions are split) and so condition (E) of $(5.4)$ is
satisfied.  Hence by
$(5.4)$, all codimension $2$ ACM schemes are glicci.

\bigskip
\noindent
{\bf Theorem 6.2.} {\em If $X$ is a nonsingular quadric three-fold in ${\mathbb
P}^4$, then two curves are in the same even $G$-liaison class if and only if
their Rao modules are isomorphic up to twist.}

\bigskip
\noindent
{\em Proof.} Indeed, for the same reason as above, condition (E$'$) of $(5.5)$ is
satisfied, and the result follows from $(5.5)$.

\bigskip
\noindent
{\bf Remark 6.3.}  One of the big open questions in the theory of Gorenstein
liaison is to describe the structure of an even $G$-liaison class
\cite[\S5.4]{M}.  Is there anything analogous to the Lazarsfeld--Rao property for
codimension $2$ CI-biliaison classes?  In this context, we can say something
about curves on the nonsingular quadric hypersurface $X$ in ${\mathbb P}^4$.  As
we have seen $(6.2)$, an even $G$-liaison class of curves on $X$ is determined
by its Rao module.  Each such even $G$-liaison class will be an infinite union of
CI-biliaison classes, each one of which has minimal curves, and satisfies the
LR-property with respect to complete intersection biliaison \cite[2.4]{RTLRP}. 
Since $\mbox{Pic } X = {\mathbb Z}$, every Gorenstein biliaison is already a
CI-biliaison, so there is no further connection among these curves using
$G$-biliaison.  However, one could hope to describe all the minimal curves in
each of these CI-biliaison classes and how they are related by $G$-liaison to
give some more structure to the whole $G$-liaison class.  All this discussion
refers to liaisons and biliaisons on $X$.  It is still possible that if one
regards curves in different CI-liaison classes on $X$, they may be related by
Gorenstein biliaisons in ${\mathbb P}^4$.

\bigskip
\noindent
{\bf Remark 6.4.}  Looking at those curves in ${\mathbb P}^4$ with Rao module $k$
that lie on the nonsingular quadric hypersurface $X$ \cite[\S4]{SEGL}, we see
that there are two types having Rao module in degree $0$, namely two skew lines,
with degree and genus $(d,g) = (2,-1)$ and the union of a conic with a line
$(3,-1)$ \cite[4.1]{SEGL}.  With Rao module $k$ in degree $1$, the $(5,0)$ curve
will  be minimal in its CI-biliaison class on $X$, while the $(6,1)$ and $(7,2)$
curves are obtained by ascending biliaison from the $(2,-1)$ and $(3,-1)$ curves,
respectively.  Note also the $(10,6)$ curve constructed at the end of
\cite[4.6]{SEGL} is not obtained by any ascending biliaison, so it is minimal in
its CI-biliaison class, with Rao module in degree $2$.  Thus we may expect to
find curves that are minimal for CI-biliaison with Rao module $k$ in arbitrarily
high degree.

\section{AG schemes in ${\mathbb P}^n$ are glicci}
\label{sec7}

\noindent
{\bf Theorem 7.1.} {\em Let $C$ be an {\em AG} subscheme of ${\mathbb P}^n$ of
codimension $\ge 2$.  Then $C$ is in the Gorenstein liaison class of a complete
intersection in ${\mathbb P}^n$.}

\bigskip
\noindent
{\em Proof.}  According to the theorem of Kleiman and Altman \cite[(1)]{KA} there
exists a complete intersection scheme $X \subseteq {\mathbb P}^n$, of dimension
$\dim C + 2$ containing $C$, that is nonsingular outside of $C$.  In particular,
$X$ will be AG, since it is a complete intersection, and normal, since it is
nonsingular in codimension $1$.

Since $C$ is an AG subscheme of $X$ of codimension $2$, it has an ${\mathcal
N}$-type resolution
\[
0 \rightarrow {\mathcal O}_X(-a) \rightarrow {\mathcal E} \rightarrow {\mathcal
I}_C \rightarrow 0
\]
where ${\mathcal E}$ is an ACM sheaf of rank $2$ on $X$ $(2.9)$.  Let ${\mathcal
M}$ be any rank $2$ dissoci\'e sheaf on $X$.  Then we get another ${\mathcal
N}$-type resolution for $C$,
\[
0 \rightarrow {\mathcal O}_X(-a) \oplus {\mathcal M} \rightarrow {\mathcal E}
\oplus {\mathcal M} \rightarrow {\mathcal I}_C \rightarrow 0.
\]
We apply $(3.4)$ to this resolution with the exact sequence
\[
0 \rightarrow {\mathcal E} \rightarrow {\mathcal E} \oplus {\mathcal M}
\rightarrow {\mathcal M} \rightarrow 0,
\]
and we find there is a codimension $2$ subscheme $D$ of $X$, in the same
$G$-liaison class as $C$, having an ${\mathcal N}$-type resolution
\[
0 \rightarrow {\mathcal L}' \rightarrow {\mathcal M} \rightarrow {\mathcal
I}_D(b) \rightarrow 0
\]
for some $b \in {\mathbb Z}$.  Since ${\mathcal M}$ is rank $2$ dissoci\'e, $D$
is a complete intersection in $X$.  Since $X$ itself is a complete intersection
in ${\mathbb P}^n$, $D$ is also a complete intersection in ${\mathbb P}^n$.  It
follows that $C$ is glicci in $X$ and hence also glicci in ${\mathbb P}^n$.

\bigskip
\noindent
{\bf Remark 7.2.} If $C$ has codimension $2$ in ${\mathbb P}^n$ then it is
already a complete intersection, so there is nothing to prove.  If $C$ has
codimension $3$, then it follows from Watanabe's paper \cite{W} that $C$ is in
fact licci (in the CI-liaison class of a complete intersection), a stronger
result.  In the special case of curves in ${\mathbb P}^4$, one can show that a
general AG curve in ${\mathbb P}^4$ can be obtained by ascending CI-biliaisons
from a line, an even stronger result \cite{GAG}.

Thus our result $(7.1)$ is new only for codimension $C \ge 4$.

\end{document}